\documentclass[12pt]{article}
\usepackage{makeidx}
\usepackage{graphicx}
\usepackage{amssymb,amsfonts,amsmath}
\usepackage{verbatim}
\usepackage[bookmarksnumbered,colorlinks,plainpages]{hyperref}

\begin{document}

\title{Composition of Fractals}

\author{Y. Lanoisel\'ee$^{1,2,3}$, L. Nivanen$^3$,\protect\\ A. El Kaabouchi$^3$ and Q. A. Wang$^{3,4}$ \thanks{Corresponding author, alexandre.wang@univ-lemans.fr} \\
{\small $^{1}$Laboratoire de Physique de la Mati\'ere Condens\'e,} \\
{\small Ecole Polytechnique, CNRS, 91128 Palaiseau Cedex, France.}\\
{\small $^{2}$Polytech Tours, 64 Avenue Jean Portalis, 37200 Tours, France.}\\
{\small $^3$Laboratoire de Physique Statistique et Syst\`emes Complexes,} \\
{\small ISMANS, 44 Avenue, F.A. Bartholdi, 72000 Le Mans, France.} \\
{\small $^4$IMMM, UMR 6283, Universit\'e du Maine}\\
{\small Avenue O. Messiaen, 72085 Le Mans Cedex, France.}}

\date{}
\maketitle{}

\begin{abstract}
This work is an analytical and numerical study of the composition of several fractals into one and of the relation between the composite dimension and the dimensions of the component fractals. In the case of composition of standard IFS with segments of equal size, the composite dimension can be expressed as a function of the component dimensions. But in the case of the compositions including component multifractals, the composite dimension cannot be expressed as explicit function of component dimensions and can only be solved numerically. An application of fractal composition to a physics problem within the incomplete statistics is discussed.
\end{abstract}

%\tableofcontents
\section{Introduction}
A fractal is a geometrical set that typically displays self-similar patterns which may be exactly or nearly the same at different scales \cite{Mand,Mand2,Falconer,Nivanen1}. In general, the concept of fractal extends beyond self-similarity and includes the case of patterns that display details at every scale. A fractal scales differently from the ordinary geometrical figures in that its box dimension $\alpha$ \cite{Falconer} is not necessarily integer. For instance, an ordinary line has $\alpha=1$, and an ordinary area has $\alpha=2$. But for a fractal line $\alpha$ can be more or less than unity. The fractal dimension is not the unique parameter describing a fractal structure, but it is the most important one characterizing the scaling property of fractals and has played important role in the application of fractal geometry to other scientific and technological domains as a modeling and design tool \cite{Mand,Mand2,Falconer}.

In this work, we focus on the mathematical operation of composing different fractals or multifractals and the calculation of the fractal dimension. By composition of fractals, we mean an operation to mix two or more fractals (multifractals) in some manner to construct a single composite fractal or multifractal. The mixing of different fractals can be performed in different manners by alternating the iterations of the component fractals in a deterministic or probabilistic way. A probabilistic mixing of different fractals has been recently proposed by Barnsley et al in the framework of the V-fractals and superfractals \cite{Barnsley1,Barnsley2} in order to facilitate the modeling of natural phenomena by superfractals. The present work is focused on the deterministic mixing of fractals. The main aim is to find the relationship between the dimensions of the component fractals and the dimension of the composite fractal, and to see whether or not this relationship is dependent on the mixing manner. This question has been raised several years ago within a statistical theory called incomplete statistics (IS) which has been proposed by physical consideration in order to make statistics in, amount others, fractal phase spaces of dynamic systems \cite{Wang2,Wang3,Wang4,Wang5,Wang6,Wang1}. A problem in this framework was to find the characteristic parameter (associated to the fractal dimension, see section 5.2 below) of a composite system from the same parameters of the component systems with the condition of the thermodynamic equilibrium \cite{Wang1,Li,Nivanen}. The second aim of this work is to justify a hypothesis of the IS (see section 5.2 below) under the angle of the composition of fractals.

\section{Composition of fractals}
A fractal can be obtained from an Iterated Function System (IFS) \cite{Falconer,Barnsley3}. An IFS is a mapping in which, starting with a closed subspace of $\mathbb{R}^n$, we iterate infinitely a set of contractions on it. A fractal is generated as the attractor of an IFS in which all contractions have the same scale factor. So, it can be described by two numbers :
 \begin{itemize}
 \item $N$ : the number of  copies (contractions) created at each stage from one copy of the previous stage with $N\in\mathbb{N}^*$.
 \item $\rho$ : the scale factor applied on each copy. We have $0<\rho< 1$ verified in the following for any IFS case and any generalization.
 \end{itemize}

\subsection{Box dimension}
\begin{figure}[!ht]
\centering
\includegraphics[width=0.5\textwidth]{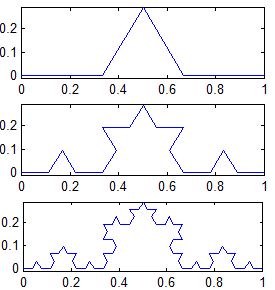}
\caption{Koch Curve ($N=4$, $\rho=1/3$) from stage $k=1$ to 3. The dimension of this curve is $\alpha=\frac{\ln 4}{\ln 3}\approx1.26$.}
\label{Koch}
\end{figure}

We start with the calculation of the box dimension for a fractal : the Koch curve (Fig.\ref{Koch}). The construction of this fractal results in an infinite curve length. In what follows, we will use the content $C_{k,\alpha}$ of the Koch curve from satge $k$. The number $C_{k,\alpha}$ is defined by the $\alpha-$ dimensional Hausdorff $\cal H^{\alpha}$ of $F_k$, where $F_k$ is the curve at the $k^{th}$ iteration \cite{Falconer}. For an ordinary line, $\alpha=1$ and for an ordinary area, $\alpha=2$. For a fractal shape, $\alpha$ is not necessarily an integer.

At each stage of the iteration of an IFS, $N$ copies are produced with a scale factor $\rho$ from every previous copy. We begin with an initial segment of length $L_0$ which is called the initiator. $F_k$ is consisted of $N_\eta=N^k$ copies of length $\eta=\rho^kL_0$.

The content $C_{k,\beta}$ is obtained with : $C_{k,\beta}=N_\eta\eta^\beta$ \cite{Falconer}. We can write : $C_{k,\beta}=(N\rho^\beta)^kL_0^{\ \beta}$. A fractal is only defined when $k\rightarrow\infty$. To reach a finite content and obtain the real dimension $\alpha$, it's necessary that $\alpha$ verifies $N\rho^\alpha=1$ which can be written as
$$\ln N+\alpha\ln\rho=0$$
or
\begin{equation} \label{Dim1IFS}
\alpha=\frac{\ln N}{\ln \frac{1}{\rho}}
\end{equation}
which is the box dimension. In the case of Fig.\ref{Koch}, $\alpha=\frac{\ln 4}{\ln 3}\approx1.26$

\subsection{Composition of 2 fractals}

The set of contractions of an IFS which is applied at each iteration is a geometric transformation. When the scale factors of the set of contractions are equals, the transformation can be expressed as a sequence of angles. Koch curve can be written as $\left[0, \frac{\pi}{3}, -\frac{\pi}{3}, 0\right]$. In what follows, we call it $K_{\frac{\pi}{3}}$. A modified Koch curve with angle $\theta$ can be written as $K_{\theta}$. For an IFS, the same transformation is applied until infinity. Koch curve can be represented by a infinite sequence of transformations as follows:

$$Koch\:Curve= K_{\frac{\pi}{3}}\:K_{\frac{\pi}{3}}\:K_{\frac{\pi}{3}}\:K_{\frac{\pi}{3}}\:K_{\frac{\pi}{3}}\:K_{\frac{\pi}{3}}\:K_{\frac{\pi}{3}}...=\underline{K_{\frac{\pi}{3}}}$$

\begin{figure}[!ht]
\centering
\includegraphics[width=0.5\textwidth]{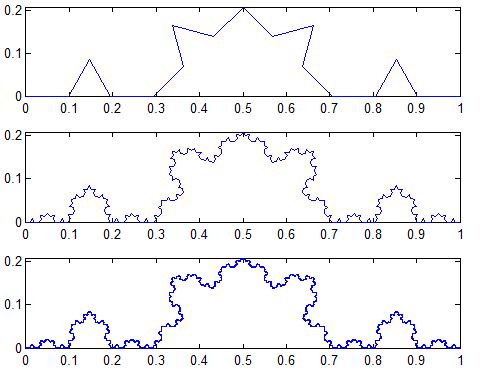}
\caption{A modified Koch curve with angles of $\frac{\pi}{4}$ ($N_1=4$, $\rho_1=\frac{1}{2+\sqrt{2}}$) composed with a Koch Curve ($N_2=4$, $\rho_2=1/3$). It's the result for $\underline{K_{\frac{\pi}{4}}K_{\frac{\pi}{3}}}$ for stage $k=1$ to 3. The component dimension of the curve is $\alpha\approx1.19$.}
\label{Kochpi4pi3}
\end{figure}

Here we are interested in composing fractals in the following way: at each stage we apply once the transformation associated with the first fractal and once the transformation associated with the second fractal.
If we want to compose Koch curve by a modified Koch curve with angle of $\frac{\pi}{4}$, the composite fractal can be represented by

$$Koch\:by\:Modified\:Koch=K_{\frac{\pi}{4}}K_{\frac{\pi}{3}}\:K_{\frac{\pi}{4}}K_{\frac{\pi}{3}}\:K_{\frac{\pi}{4}}K_{\frac{\pi}{3}}...=\underline{K_{\frac{\pi}{4}}K_{\frac{\pi}{3}}}$$
where each transformation is called a substage and a period of transformation is called a stage. Hence each stage is composed of a $K_{\frac{\pi}{4}}$ substage and a $K_{\frac{\pi}{3}}$ substage. The underline means the couple of transformations to be repeated until infinity. Each stage associated with the composite fractal produces respectively $N_1$ and $N_2$ copies with respectively scale factors $\rho_1$ and $\rho_2$ (Fig.\ref{Kochpi4pi3}).

Hence at the first substage, $N_1$ copies are generated with a scale factor $\rho_1$. And at the second substage $N_2$ copies are generated with a scale factor $\rho_2$ from  each of the $N_1$ copies of the first substage. So at stage $k$ we have $N_\eta=N_1^{\ k}N_2^{\ k}$ copies of length $\eta=\rho_1^{\ k}\rho_2^{\ k}L_0$. The content is
\begin{equation} \label{Content2IFS}
C_{k,\beta}= (N_1N_2(\rho_1\rho_2)^\beta)^kL_0^{\ \beta}
\end{equation}
To have a finite content when $k\rightarrow\infty$, it's necessary that :
$N_1N_2(\rho_1\rho_2)^\alpha=1$, which leads to :

\begin{equation} \label{Dim2IFS}
\alpha=\frac{\ln N_1 + \ln N_2}{\ln \frac{1}{\rho_1}+\ln \frac{1}{\rho_2}}
\end{equation}

The dimension of the composite fractal can be expressed as a function of the dimensions of the component fractals. Remember that each component fractal verifies the relation $N_i\rho_i^{\alpha_i}=1$ with $i\in\{1,2\}$ where $\alpha_i$ is its dimension, we can write $N_i=\frac{1}{\rho_i^{\alpha_i}}$ and replace it in the above equation to give
 $\alpha=\frac{\ln{\frac{1}{\rho_1^{\alpha_1}}} + \ln{\frac{1}{\rho_2^{\alpha_2}}} }{\ln \frac{1}{\rho_1}+\ln \frac{1}{\rho_2}}$ and :
 \begin{equation} \label{Dim(Dim)2IFS}
 \alpha=\frac{ \alpha_1 \ln \rho_1 +  \alpha_2 \ln \rho_2}{\ln \rho_1+\ln \rho_2}
\end{equation}

The composite dimension is independent of the order of the substages in the IFS. Note that its expression is in the form of a barycentric average with positive coefficients $\frac{\ln \rho_1}{\ln \rho_1+\ln \rho_2}$ and $\frac{\ln \rho_2}{\ln \rho_1+\ln \rho_2}$. Because of this average, we know that $\alpha$ lies between  $\alpha_1$ and $\alpha_2$. For Fig.\ref{Kochpi4pi3}, using Eq.\eqref{Dim2IFS}, we can calculate $\alpha\approx1.19$.

\subsection{Composition of \texorpdfstring{$m$}{TEXT} fractals}
This section describes how to obtain the dimension of a composite fractal composed of $m$ different fractals. We construct it by an IFS of $m$ substages (corresponding to $m$ sets of contractions) in each stage of the iteration. Let $T_i$ be the transformation associated with the $i^{th}$ fractal, we can represent the composition with:
$$Composition\:of\:m\:fractals=T_1T_2...T_mT_1T_2...T_mT_1T_2...T_m...=\underline{T_1T_2...T_m} $$

\subsubsection{General Case}

\begin{figure}[!ht]
\centering
\includegraphics[width=0.5\textwidth]{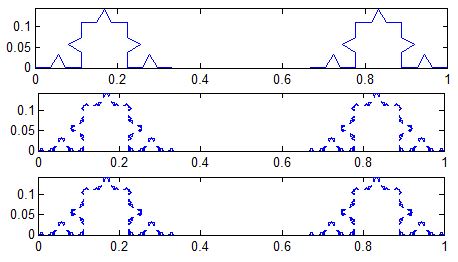}
\caption{A composite IFS composed of a Cantor set ($N_1=2$, $\rho_1=1/3$), a quadratic Koch curve ($N_2=5$, $\rho_2=1/3$) and a Koch curve ($N_3=4$, $\rho_3=1/3$). The composition rule is  $Composition=\underline{C_{[1/3\ 1/3]}Q_{\frac{\pi}{2}}K_{\frac{\pi}{3}}}$ for the stage $k$ with $k=1$ to 3. The composite dimension is $\alpha\approx1.12$.}
\label{Cantorpi2pi3}
\end{figure}

We make $m$ substages with $N_i$ copies of scale factor $\rho_i$ for the substage $i$ with $i\in[[1,m]]$. Each substage verifies $N_i\rho_i^{\ \alpha_i}=1$ with $\alpha_i$ the dimension of the $i^{th}$ component fractal. Fig.\ref{Cantorpi2pi3} shows an example combining a Cantor set, a quadratic Koch curve (noted $Q_{\frac{\pi}{2}}$) and a Koch curve.

At the $k^{th}$ stage, $N_\eta=\prod_{i=1}^m N_i^k$ copies of length $\eta=\left(\prod_{i=1}^m \rho_i^k\right)L_0$ are produced. The content is:

\begin{eqnarray} \label{ContentMIFS} \nonumber
 C_{k,\beta} &=&( \prod_{i=1}^m N_i\rho_i^{\ \beta})^kL_0^{\ \beta}  \\
\end{eqnarray}
When $k\rightarrow\infty$, it's necessary that :

\begin{equation} \label{RelDimMIFS}
 \prod_{i=1}^m N_i\rho_i^{\ \alpha}=1
\end{equation}
or
\begin{equation} \label{DimMIFS} \nonumber
\alpha = \frac{\sum\limits_{i=1}^m \ln N_i}{\sum\limits_{i=1}^m \ln \frac{1}{\rho_i}}
\end{equation}
For the example in Fig.\ref{Cantorpi2pi3}, from Eq.\eqref{DimMIFS} : $\alpha=\approx1.12$.

Because of the commutativity of the sum in Eq.\eqref{DimMIFS}. The dimension $\alpha$ is independent of the order of application of different transformations $T_i$, with $i\in [[1,m]]$ in a stage. The Eq.\eqref{DimMIFS} is a generalization of the result obtained in Eq.\eqref{Dim2IFS}.

To express the composite dimension as a function of the component dimensions and scale factors, we use $N_i\rho_i^{\ \alpha_i}=1$ to replace $N_i$ in Eq.\eqref{DimMIFS} and obtain

\begin{equation} \label{Dim(Dim)MIFS} \nonumber
\alpha = \frac{\sum\limits_{i=1}^m \alpha_i\ln \rho_i}{\sum\limits_{i=1}^m \ln \rho_i}
\end{equation}
If the $m$ component IFS have the same scale factor, we can write $\rho_i=\rho$.Eq.\eqref{Dim(Dim)MIFS} becomes

\begin{equation} \label{ArithmAverage}
\alpha=\frac{ \sum_{i=1}^m \alpha_i}{m} \nonumber
\end{equation}
which is an arithmetic average of the component dimensions.

We can also express the composite dimension as a function of the component dimensions and the number of copies of each substage $N_i$. Let us take Eq.\eqref{DimMIFS} and replace $\rho_i$ by $\rho_i=\frac{1}{N_i^{\ \frac{1}{\alpha_i}}}$. The composite dimension then verifies

\begin{eqnarray} \label{Dim2(Dim)MIFS}
\frac{1}{\alpha} &=&  \frac{\sum\limits_{i=1}^m \frac{1}{\alpha_i}\ln N_i}{\sum\limits_{i=1}^m \ln N_i} \\ \nonumber
\end{eqnarray}
where the inverse composite dimension is in the form of a barycentric average of the component dimension determined by the integer numbers $N_i$.

If the $m$ component fractals have the same number of copies, we can write $N_i=N$ and obtain
$$\frac{1}{\alpha} = \frac{\sum\limits_{i=1}^m \frac{1}{\alpha_i}}{m} $$
i.e., the composed dimension is the harmonic average of the component dimensions.

\section{Composition of multifractals}
A multifractal set is the attractor of an IFS where the scale factors can be different. At a given stage, the segments have different sizes and forms which can be treated with multinomial expansions. Let $r_j$ be the multifractal scale factors with $j\in [[1,l]]$ where $l$ is the number of scaling applied. As in the case of fractals, we shall verify $\forall j \in [[1,l]]$, $0<r_j< 1$ for any multifractal. To have non zero finite content, each multifractal verifies the Moran equation\cite{Wang1}:
$$\sum\limits_{j=1}^lr_j^{\ \alpha}=1$$
with $\alpha$ the dimension of the multifractal. If all scale factors are equal $r_j=r$ , the Moran equation becomes $lr^\alpha=1$, the relation for standard IFS.

It is possible to prove the uniqueness of $\alpha$. We take a function $f(\alpha)=\sum\limits_{j=1}^lr_j^{\ \alpha}$. It is a continuous and constantly decreasing function. When $\alpha\rightarrow 0$, $f(\alpha)\rightarrow l$ and when $\alpha\rightarrow\infty$, $f(\alpha)\rightarrow 0$. So $\sum\limits_{j=1}^lr_j^{\ \alpha}=1$ has unique solution corresponding to the box dimension of the curve.

\subsection{Composition of a fractal and a binary multifractal}

\begin{figure}[!ht]
\centering
\includegraphics[width=0.5\textwidth]{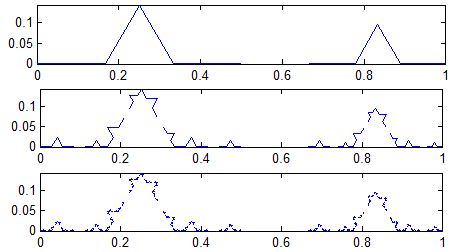}
\caption{Composition of a binary multifractal ($r_1=1/2$, $r_2=1/3$) with a Koch Curve ($N=4$, $\rho=1/3$) The composition rule is $Composition = \underline{C_{[1/2\ 1/3]}K_\frac{\pi}{3}}$ for the stage $k$ with $k=1$ to 3.}
\label{CantorAsymKoch}
\end{figure}

A binary multifractal is constructed with two different scale factors $r_1$ and $r_2$, with $\{r_1,r_2\}<1$ and $r_1^{\ \alpha_r}+r_2^{\ \alpha_r}=1$ \cite{Wang1}, where $\alpha_r$ is the dimension of the binary multifractal. In what follows, we focus on the composition of a binary multifractal with a fractal. Fig.\ref{CantorAsymKoch} gives an example of composition of a binary multifractal and a Koch curve.

Each stage of the iteration is composed of two substages: one being the iteration of the binary multifractal and another one is the iteration of the fractal. At the $k^{th}$ stage $F_k$ is composed of $N_k=2^kN^k$ copies. At a given scale, the length of a segment is $\delta_{k,p}$ with $p\in [[1,N_k]]$, so the total length of the curve is $L_k=\sum\limits_{p=1}^{N_k}\delta_{k,p}$. Now consider $j\in [[0,k]]$,
we can find $ \binom{k}{j} $ segments of length $\Delta_{k,j}=r_1^{\ j} r_2^{\ k-j} L_{0}$ on the binary multifractal\cite{Wang1}. Each segment of this length contains $N^k$ segments of length $\rho^kL_0$ generated by fractal. It is obvious that at the $k^{th}$ stage we have $\binom{k}{j} N^k$ segments of length $\Delta_{k,j}=r_1^{\ j}r_2^{\ k-j}\rho^kL_0$. The total length of the curve is then:
\begin{equation}
L_k=\sum\limits_{p=1}^{N_k}\delta_{k,p}=\sum_{j=0}^k\binom{k}{j}N^kr_1^{\ j}r_2^{\ k-j}\rho^kL_0=N^k \rho^k \left(r_1+r_2\right)^kL_0\nonumber
\end{equation}
The content of this form is given by:
\begin{equation}
C_{k,\beta}=\sum\limits_{p=1}^{N_k}\delta_{k,p}^{\ \ \ \beta}=\sum\limits_{j=0}^k\binom{k}{j}r_1^{\ \beta j}r_2^{\ \beta(k-j)}N^k\rho^{\beta k}L_0^{\ \beta}
\end{equation}
in which the binomial expression is
$$\sum\limits_{j=0}^k\binom{k}{j}r_1^{\ \beta j}r_2^{\ \beta(k-j)} =(r_1^{\ \beta}+r_2^{\ \beta})^k $$
Hence the content of order $\beta$ reads:
\begin{equation}
C_{k,\beta}=\left[(r_1^{\ \beta}+r_2^{\ \beta})N\rho^\beta\right]^kL_0^{\ \beta}
\end{equation}
The content must be finite when $k\rightarrow\infty$, implying
\begin{equation}
(r_1^{\ \alpha}+r_2^{\ \alpha})N\rho^\alpha=1
\end{equation}
This expression gives a unique $\alpha$ as the composite dimension of the curve. The dependence of $\alpha$ on the component dimension $\alpha_\rho$ of the fractal can be shown as follows. $\alpha_\rho$ verifying $N\rho^{\alpha_\rho}=1 \Leftrightarrow N=\frac{1}{\rho^{\alpha_\rho}}$, the above equation can be written as
\begin{equation}\label{Dim(Dim)BinaryIFS}
(r_1^{\ \alpha}+r_2^{\ \alpha})\rho^{\alpha-\alpha_{\rho}}=1
\end{equation}
However, in general it is impossible to introduce the multifractal dimension $\alpha_r$. $\alpha$ should be obtained by numerical solution of the above two equations.

Anyway, analytical solution of $\alpha$ from Eq.\eqref{Dim(Dim)BinaryIFS} is possible under special conditions. For example, one of these conditions is $r_2=r_1^2\rho$. The proof is the following.

Let us put $r_2$ into Eq.\eqref{Dim(Dim)BinaryIFS} to get

   \[(r_1^{\ \alpha}+r_1^{\ 2\alpha}\rho^\alpha)\rho^{\alpha-\alpha_\rho}=1 \] 
or
   \[(r_1^{\ \alpha}\rho^\alpha)^2+r_1^{\ \alpha}\rho^\alpha-\rho^{\alpha_\rho}=0\] 

Because $(r_1\rho)^\alpha>0$ and $\rho^{\alpha_\rho}=\frac{1}{N}$ it follows that
   \[ (r_1\rho)^\alpha=\frac{-1+\sqrt{1+\frac{4}{N}}}{2}\]

   Finally:

 \begin{equation} \label{DimBinAnal}
 \alpha=\frac{\ln\frac{-1+\sqrt{1+\frac{4}{N}}}{2}}{\ln r_1\rho}
 \end{equation}

\begin{figure}[!ht]
\centering
\includegraphics[width=0.5\textwidth]{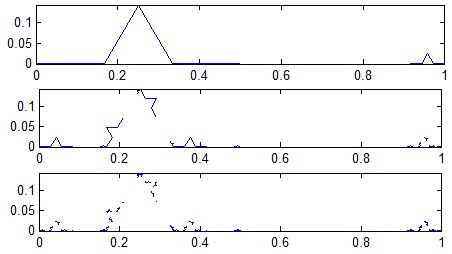}
\caption{A binary multifractal ($r_1=1/2$, $r_2=1/12$) composed with Koch curve ($N=4$, $\rho=1/3$, $\alpha_\rho\approx1.26$). It's the result for $\underline{C_{[1/2\ 1/12]}K_{\frac{\pi}{3}}}$ for stage $k=1$ to 3. The component dimension of the curve is $\alpha\approx0.88$.}
\label{CantorAsymAnalytic}
\end{figure}

Fig.\ref{CantorAsymAnalytic} shows an example of composition of a Cantor ($r_1=1/2$, $r_2=1/12$) set with a Koch curve ($N=4$, $\rho=1/3$).

where $r_2=\left(\frac{1}{2}\right)^2\frac{1}{3}=\frac{1}{12}$.\\
 The dimension is given by ${\alpha=\frac{\ln\frac{-1+\sqrt{1+\frac{4}{4}}}{2}}{\ln \left(\frac{1}{2}\frac{1}{3}\right)}\approx 0.88}$.

\subsection{Composition of \texorpdfstring{$m$}{TEXT} fractals and a multifractal}

\begin{figure}[!ht]
\centering
\includegraphics[width=0.5\textwidth]{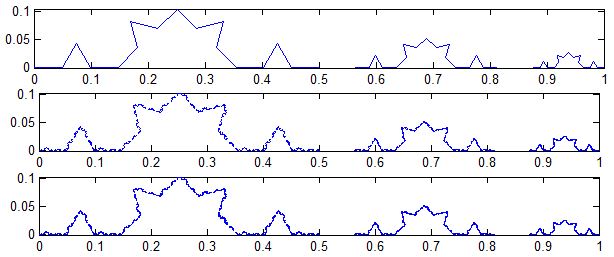}
\caption{Composition, from stage 1 to 3, of a multifractal (${r_1=1/2}$, $r_2=1/4$, $r_3=1/6$), a modified Koch curve for $\theta=\pi/4$ ($N_1=4$, $\rho_1=\frac{1}{2+\sqrt{2}}$) and a Koch curve (${N_2=4}$, ${\rho_2=1/3}$). This composition is represented by :\newline
${Composition=C_{\left[1/2,1/4,1/6\right]}K_\frac{\pi}{4}K_\frac{\pi}{3}C_{\left[1/2,1/4,1/6\right]}K_\frac{\pi}{4}K_\frac{\pi}{3}...=\underline{C_{\left[1/2,1/4,1/6\right]}K_\frac{\pi}{4}K_\frac{\pi}{3}}}$.}
\label{MultiCantorKochpi3Kochpi4}
\end{figure}

Now we focus on the composition of $m$ different fractals and a multifractal. At the $k^{th}$ stage, each one of the $m$ transformations associated with fractals produce $N_i$ copies of scale factor $\rho_i$ with $i\in [[1,m]]$, with $ \prod\limits_{i=1}^m N_i^{\ k} $ the total number of segments of scale factor $ \prod\limits_{i=1}^m \rho_i^{\ k}$. On the other hand, the multifractal set is defined by $l$ different scale factors $r_j$ with $r_j<1$ and $j\in [[1,l]]$. An example of this composition is given in Fig.\ref{MultiCantorKochpi3Kochpi4}.

At the $k^{th}$ stage, on $F_k$ we find $N_k=l^k\prod\limits_{i=1}^mN_i^{\ k}$ segments of length $\delta_{k,p}$ with $p\in[[1,N_k]]$. Let $G$ be any set of integers $\{g_1,g_2,...,g_m\}$ verifying :
\begin{equation}\label{multinomeG}
 \forall j\in [[1,l]],\ g_j\in [[0,k]],\  \sum_{j=1}^l g_j = k
\end{equation}

For any set G that verifies  Eq.\eqref{multinomeG}, we can find, on a multifractal, $ \frac{k!}{\prod\limits_{j=1}^l g_j !}$ segments of length $ \Delta_{k,G}=(\prod\limits_{j=1}^l r_j^{\ g_j})L_0$ \cite{Wang1}

At the $k^{th}$ stage, on the $F_k$ form, we have $ \frac{k!}{\prod\limits_{j=1}^l g_j !}\prod\limits_{i=1}^m N_i^{\ k}$ segments of length : 
$ \Delta_{k,G}=\left(\prod\limits_{j=1}^l r_j^{\ g_j}\right) \left( \prod\limits_{i=1}^m \rho_i^{\ k}\right) L_0 $.

Then the total length is :

\begin{align}
\nonumber
 L_k &= \sum_{p=1}^{N_k} \delta_{k,p} \\ \nonumber
 &= \sum_{G}\left(\frac{k!}{\prod\limits_{j=1}^l g_j !}\right)\left(\prod_{j=1}^l r_j^{\ g_j}\right)\left(\prod_{i=1}^m N_i \rho_i\right)^kL_0 \\ 
L_k &= \left(\sum\limits_{j=1}^lr_j\right)^k\left(\prod\limits_{i=1}^m N_i \rho_i\right)^kL_0
\end{align}

Now let us calculate the content of order $\alpha$ on $F_k$ :
$$ C_{k,\beta}=\sum\limits_{p=1}^{N_k} \delta_{k,p}^{\ \ \ \beta}= \sum\limits_{G}\left(\frac{k!}{\prod\limits_{j=1}^l g_j !}\right)\left(\prod\limits_{j=1}^l r_j^{\ \beta g_j}\right)\left(\prod\limits_{i=1}^m N_i \rho_i^{\ \beta}\right)^kL_0^{\ \beta} $$
We have here a multinomial expansion.
\begin{equation}
 \sum\limits_{G}\left(\frac{k!}{\prod\limits_{j=1}^l g_j !}\right)\left(\prod\limits_{j=1}^l r_j^{\ \beta g_j}\right) = \left(\sum_{j=1}^l r_j^{\ \beta} \right)^k
 \end{equation}
It leads to :
\begin{equation}\label{ContentMultiMIFS}
C_{k,\beta}=\left[ \left(\prod_{i=1}^m N_i\rho_i^{\ \beta}\right) \left( \sum_{j=1}^l r_j^{\ \beta} \right) \right]^k L_0^{\ \beta}
\end{equation}
When $k\rightarrow\infty$ the content remains finite with the dimension verifying :
\begin{equation}\label{DimMultiMIFS}
\left(\prod_{i=1}^m N_i\rho_i^{\ \alpha} \right) \left( \sum_{j=1}^l r_j^{\ \alpha} \right)=1
\end{equation}
or
\begin{equation}\label{DimMultiMIFS2}
\left(\prod_{i=1}^m \rho_i^{\ \alpha-\alpha_{\rho,i}} \right) \left( \sum_{j=1}^l r_j^{\ \alpha} \right)=1.
\end{equation}

\subsection{Composition of \texorpdfstring{$m$}{TEXT} multifractals}

\subsubsection{Composite dimension}
In this section we study the composition of $m$ different multifractals.
Each one of these multifractals is defined by $l_i$ scale factors with $i\in [[1,m]]$ corresponding to the sets of scale factors from each one of the $m$ multifractals.

At the $k^{th}$ stage, $F_k$ is composed of $N_k=\prod\limits_{i=1}^m l_i^{\ k}$ segments of length $\delta_{k,p}$ with $p\in [[1,N_k]]$.

We note $G_i$ any set of integers $\{g_{i,1},g_{i,2},\cdots,g_{i,l_i}\}$ that verifies:
\begin{equation} \label{multinomeG1}
 \forall j\in \{1\cdots l_i\},\  g_{i,j}\in\{0,1,\cdots,k\},\ \sum_{j=1}^{l_i} g_{i,j} = k
\end{equation}

For any set $G_i$ that verifies Eq.\eqref{multinomeG1} : \\
We can find $ \frac{k!}{\prod\limits_{j=1}^{l_i} g_{i,j} !}$ segments of length $ \Delta_{k,G_i}=(\prod\limits_{j=1}^{l_i} r_{i,j}^{\ \ \ g_{i,j}})L_0$.

The total length of the $F_{k}$ form is given by:
$$L_k=\sum_{p=1}^{N_k} \delta_{k,p} =\prod\limits_{i=1}^m\left[\sum_{G_i}\left(\frac{k!}{\prod\limits_{j=1}^{l_i} g_{i,j} !}\right)\left(\prod_{j=1}^{l_i} r_{i,j}^{\ \ \ g_{i,j}}\right)\right]L_0=\prod\limits_{i=1}^m\left(\sum\limits_{j=1}^{l_i} r_{i,j}\right)^kL_0$$

The content of order $\beta$ of the form $F_k$ is:
\[C_{k,\beta}=\sum\limits_{p=1}^{N_k} \delta_{k,p}^{\ \ \ \beta}= \prod\limits_{i=1}^m\left[\sum\limits_{G_i}\left(\frac{k!}{\prod\limits_{j=1}^{l_i} g_{i,j} !}\right)\left(\prod\limits_{j=1}^{l_i} r_{i,j}^{\ \ \ \beta g_{i,j}}\right)\right]L_0^{\ \beta} \]

This is the form of a product of multinomial expansions.
$$\prod\limits_{i=1}^m\left[\sum\limits_{G_i}\left(\frac{k!}{\prod\limits_{j=1}^{l_i} g_{i,j} !}\right)\left(\prod\limits_{j=1}^{l_i} r_{i,j}^{\ \ \ \beta g_{i,j}}\right)\right] = \prod\limits_{i=1}^m\left(\sum\limits_{j=1}^{l_i} r_{i,j}^{\ \ \ \beta}\right)^k $$
The content is recast as:
\[C_{k,\beta}= \prod\limits_{i=1}^m\left(\sum\limits_{j=1}^{l_i} r_{i,j}^{\ \ \ \beta}\right)^kL_0^{\ \beta} \]

Because the content is finite when $k\rightarrow\infty$, we have
\begin{equation} \label{Mmultifractal}
\prod\limits_{i=1}^m\left(\sum\limits_{j=1}^{l_i} r_{i,j}^{\ \ \ \alpha}\right)=1
\end{equation}
As all $r_{i,j}<1$, the expression is strictly decreasing so there is a unique $\alpha$ which satisfies Eq.\eqref{Mmultifractal}.
Since multifractal is a general case of fractals, Eq.\eqref{Mmultifractal} is a general formula of composite dimension for all compositions.

\subsubsection{Borders of the composite dimension}

Since the composition of $m$ multifractals is a generalization of the composition of $m$ fractals, it is reasonable to think of the existence of borders for the composite dimension such as $\underset{i\in[[1,m]]}{min(\alpha_i)}\leq \alpha \leq \underset{i\in[[1,m]]}{max(\alpha_i)}$. This can be proved by contradiction.
Let's postulate $\underset{i\in[[1,m]]}{max(\alpha_i)}<\alpha$.
The $m$ multifractals verify $\sum\limits_{j=1}^{l_i}r_{i,j}^{\ \ \ \alpha_i}=1$ with $i\in [[1, m]]$.
Let $f_i(\alpha_i)$ be a function of the component dimensions $\alpha_i$ such that $f_i(\alpha_i)=\sum\limits_{j=1}^{l_i}r_{i,j}^{\ \ \ \alpha_i}$. Then all $f_i$ are strictly decreasing functions so:
$$\forall i\in [[1,m]],\ f_i(\alpha )<f_i(\alpha_i)=1$$

Hence $$\prod\limits_{i=1}^mf_i(\alpha)<\prod\limits_{i=1}^mf_i(\alpha_i) =1$$
But from Eq.\eqref{Mmultifractal}, we know that $\prod\limits_{i=1}^mf_i(\alpha)=1$.
This is a contradiction, hence we must have

\begin{equation} \label{Encadrement_multifractal}
\underset{i\in[[1,m]]}{min(\alpha_i)}\leq\alpha\leq\underset{i\in[[1,m]]}{max(\alpha_i)}
\end{equation}
We can conclude that the composite dimension is always between the minimum and the maximum of component dimensions for any composition.

\section{Periodically alternating compositions}
The results presented above concerns the compositions performed by alternating the iterations of component fractals one after another, meaning that in each stage of the composite iteration, the number of substages is equal to the number of different component fractals. In what follows, we will study a more general case where, in a stage of the composite iteration, each component fractal is repeated several times following by repeating the iteration of another component fractals.

\subsection{For 2 fractals}
Consider two component fractals whose corresponding transformations produces in each iteration $N_1$ and $N_2$ copies of scale factor $\rho_1$ and $\rho_2$, respectively. At each stage of the composition, we apply $n_1$ times the first scaling and $n_2$ times the second scaling. So each stage is composed of $n=n_1+n_2$ substages. In the following Figure, we apply once Cantor and twice Koch transformation in each stage as represented by the following sequence:
$$Cantor_{\left[1/3,1/3\right]}\:by\:Koch_{\frac{\pi}{3}}=CKK\:CKK\:CKK\:CKK\:CKK...=\underline{CKK}$$

\begin{figure}[!ht]
\centering
\includegraphics[width=0.5\textwidth]{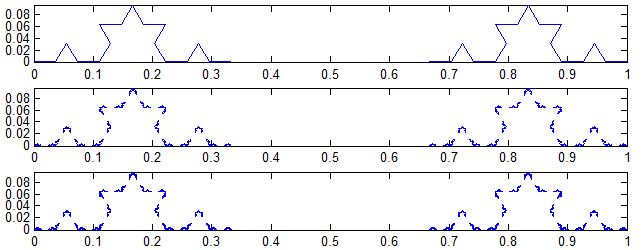}
\caption{Cantor applied once and Koch curve applied twice by stage,that makes 3 substages, with respectively ($n_1=1$, $N_1=2$, $\rho_1=1/3$) and (${n_2=2}$, $N_2=4$, $\rho_2=1/3$) from stage $k=1$ to 3. The dimension obtained from Eq.\eqref{Dim(Dim)Comp2IFS} is $\alpha\approx1.05$.}
\label{Cantor2Koch}
\end{figure}

We can use Eq.\eqref{Dim(Dim)MIFS} to find the composite dimension as if we had $n=n_1+n_2$ different substages (component fractals):
\begin{equation} \nonumber
\alpha=\frac{\sum_{i=1}^n \alpha_i\ln \rho_i}{\sum_{i=1}^n \ln \rho_i}.
\end{equation}
It is clear that in each stage of $n$ substages, there are $n_1$ times the first scaling and $n_2$ times the second one. Hence we have

\begin{equation} \label{Dim(Dim)Comp2IFS}
\alpha=\frac{n_1\alpha_1\ln \rho_1 + n_2\alpha_2\ln \rho_2}{n_1 \ln \rho_1 + n_2 \ln \rho_2}
\end{equation}
Hence we conclude that the composite dimension is independent of the order of substages in a stage if $n_1$ and $n_2$ are given. Consequently, the example in Fig.\ref{Cantor2Koch} has the same dimension as the following sequence:
$$Cantor_{\left[1/3,1/3\right]}\ by\ Koch_{\frac{\pi}{3}}=KCK\ KCK\ KCK\ KCK\ KCK\cdots=\underline{KCK}$$

Further more, the order of substages can be different from one stage to another without affecting the composite dimension if $n_1$ and $n_2$ are given. Hence the following sequence has the same composite dimension as in Fig.\ref{Cantor2Koch}.

$$Cantor_{\left[1/3,1/3\right]}\ by\ Koch_{\frac{\pi}{3}}=KCK\ CKK\ KKC\ CKK\ KCK\cdots$$

\subsection{For \texorpdfstring{$m$}{TEXT} fractals}
In this part we will extend the case of a periodic composition of 2 fractals to $m$ fractals, by applying $n_i$ times one of the $m$ transformation to produce $N_i$ segments of scale factor $\rho_i$. For example we can make a composition of 3 fractals whose iterations are represented by $T_1$, $T_2$ and $T_3$, respectively. $T_1$ is repeated $n_1=3$ times,  $T_2$ $n_2=2$ times and $T_3$ $n_3=1$ times. Hence there are 6 substages in each stage of the composition. The composition can be represented by the following sequence:
$$Composition=T_1T_1T_1T_2T_2T_3\ T_1T_1T_1T_2T_2T_3 \ldots=\underline{T_1T_1T_1T_2T_2T_3}$$
So in general each stage is composed of $\sum\limits_{i=1}^m n_i =n$ substages.
By using Eq.\eqref{Dim(Dim)MIFS}, we obtain

\begin{equation} \label{Dim(Dim)CompMIFS} \nonumber
\alpha = \frac{\sum\limits_{i=1}^n \alpha_i\ln \rho_i}{\sum\limits_{i=1}^n \ln \rho_i}
= \frac{\sum\limits_{i=1}^m n_i\alpha_i\ln \rho_i}{\sum\limits_{i=1}^m n_i \ln \rho_i }
\end{equation}

If all $n_i$ are the same, we find Eq.\eqref{DimMIFS}. It is obvious that the composite dimension, for given fractal, depends on the number of repetition of each IFS in each stage of the composition.

\subsection{For \texorpdfstring{$m$}{TEXT} multifractals}
We can also make, in the same way as above, periodic composition for multifractals. Each one of the $m$ multifractal transformation is applied $n_i$ times in each stage of the composition. The total number of substages is then $\sum\limits_{i=1}^m n_i=n$.
For example, 3 multifractals, $M_1$, $M_2$ and $M_3$, are composed periodically. Suppose ${n_1=1}$, ${n_2=2}$ and ${n_3=3}$. There are 6 substages in each stage of the composition.
The sequence is :
\begin{align}
\nonumber
Composition&=M_1M_2M_2M_3M_3M_3\ M_1M_2M_2M_3M_3M_3\ ... \\ 
&=\underline{M_1M_2M_2M_3M_3M_3}
\end{align}
The general formula for this composition is:
\begin{eqnarray} \nonumber
\prod\limits_{i=1}^m\left(\sum\limits_{j=1}^{l_i} r_{i,j}^{\ \ \ \alpha}\right)^{n_i}=1
\end{eqnarray}
Because the relation is commutative, the composite dimension does not depend on the order of the different substages. It depends only on the numbers $n_i$.

\section{Application}

\subsection{How to get a dimension close to a rational number}

Starting with an initial fractal of any dimension, we can get a composite fractal with dimension as close as we want to a rational number. Consider a fractal constructed by a transformation which creates $N_1$ copies at each iteration with a scale factor $\rho_1$. It is composed with another fractal whose transformation creates $N_2=n^{a_1}$ copies with of scale factor $\rho_2=\frac{1}{n^{a_2}}$ where $\{n,a_1,a_2\}\in\mathbb{N}^*$. The composite dimension $\alpha$ is given by

\begin{equation} \nonumber
\alpha_{lim} = \lim\limits_{n \to \infty} \frac{\ln N_1+\ln n^{a_1}}{\ln \frac{1}{\rho_1}+\ln n^{a_2}}
= \lim\limits_{n \to \infty} \frac{\frac{\ln N_1}{\ln n^{a_1}}+1}{\frac{\ln \frac{1}{\rho_1}}{\ln n^{a_1}}+ \frac{\ln n^{a_2}}{\ln n^{a_1}}}
 = \frac{a_1}{a_2}
\end{equation}
So it is possible to get a dimension as near as possible (depending on $n$) to a rational number $\frac{a_1}{a_2}$.

\subsection{Connection with Incomplete Statistics}
We study here the link between fractal composition and the incomplete statistics \cite{Wang2} which has been proposed several years ago in order to make statistical study of physical systems for which the mathematical calculation of available states is incomplete (under or over counted). One of the possible applications is to the cases where the phase space is fractal or the probability measure can be represented in fractal geometry. For instance, for a multifractal phase space, the incomplete probability distribution is defined for any segment $i_k$ at any stage $k$ by $p_{i_k}=\delta_{i_k}/L_0$ with the normalization $\sum\limits_{i=1}^lp_{i_k}^{\ \alpha}=1$ where $\alpha$ is a ratio equal to the box dimension $d_f$ divided by the dimension $d$ of the space into which the fractal is embedded: $\alpha=d_f/d$\cite{Wang1}. $\alpha$ characterizes the incompleteness of the statistics on that fractal phase space. As this incomplete normalization is valid for any $k$, we will drop the index $k$ from the normalization in what follows.

The incomplete statistics has been proposed in the context of the study of a nonextensive statistics\cite{Wang2,Wang3}. In an effort to apply the nonextensive statistics to thermodynamic systems composed of two independent subsystems, a problem is encountered about the derivation of the zeroth law of thermodynamics relative to the establishment of thermodynamic equilibrium between the two subsystems, say, $A$ and $B$, with probability distributions $p_i(A)$ and $p_j(A)$, respectively\cite{Wang7}. In order to reach the necessary zeroth law, the following hypotheses have been imposed. 1) If the probability distribution is complete and normalized with $\sum\limits_{i=1}p_i=1$ for the composite system and the subsystems, the joint probability of the composite system should be given by ${p_{ij}^\alpha(A+B)=p_{i}^{\ \alpha_a}(A)p_{j}^{\ \alpha_b}(B)}$ where ${\alpha}$, ${\alpha_a}$ and ${\alpha_b}$ are the incompleteness exponents of the composite system, subsystem $A$ and subsystem $B$, respectively\cite{Li}. 2) If the probability distribution is incomplete and normalized with $\sum\limits_{i,j}p_{ij}^{\ \ \alpha}(A+B)=1$, ${\sum\limits_{i}p_{i}^{\ \alpha_a}(A)=1}$ and $\sum\limits_{j}p_{j}^{\ \alpha_b}(B)=1$, the joint probability of the composite system should be given by ${p_{ij}(A+B)=p_{i}(A)p_{j}(B)}$\cite{Nivanen}.

In those works, the expressions of the joint probability are nevertheless hypothetical in order to get coherent results concerning the existence of the zeroth law in those thermostatistics formalisms. In what follows, we will show that the relationship of the incomplete joint probability mentioned above is a natural consequence of the composition of the fractals with independent scale factors.

Now imagine two multifractal phase spaces $A$ and $B$ in which we have the following incomplete normalizations :
\begin{equation}\label{normi}
\sum\limits_{i}p_i^{\ \alpha_a}=1
\end{equation}
and
\begin{equation}\label{normj}
\sum\limits_{j}p_j^{\ \alpha_b}=1.
\end{equation}
where the sums are over all the segments of the two fractals, respectively.

Suppose now that $A$ and $B$ compose a composite multifractal $C$ in the way described in section 3.4, we naturally have $r_{i,j}=r_{i}r_{j}$ and, at $k$ stage,  $\delta_{i_k,j_k}=\delta_{i_k}\delta_{j_k}/L_0$ or $p_{i,j}=p_ip_j$. From Eq.\eqref{Mmultifractal}, we can write
\begin{eqnarray}\label{composedSubSystems}\nonumber
\sum\limits_{i,j} p_i^{\ \alpha}p_j^{\ \alpha}=1
\end{eqnarray}
or
\begin{eqnarray}\label{compositeSystems}\nonumber
\sum\limits_{i,j} p_{ij}^{\ \ \alpha}=1
\end{eqnarray}
which is the incomplete normalization for the composite systems. Hence the hypothesis of the incomplete joint probability made in the context of nonextensive statistics\cite{Nivanen} is supported from the point of view of fractal composition presented in this work.

If a system $S$ is composed of $m$ independent subsystems $\left\{A,B,...,M\right\}$
with $\left\{l_1,l_2,...,l_m\right\}$ the indices of the probability distributions of the subsystems, from Eq.\eqref{Mmultifractal}, we can write

\begin{eqnarray} \nonumber
p_{l_1,l_2,...,l_m}\left(S\right)&=&p_{l_1}\left(A\right)p_{l_2}\left(B\right)...p_{l_m}\left(M\right)
\end{eqnarray}
and
\begin{eqnarray} \nonumber
\sum\limits_{l_1,l_2,...,l_m}\left(p_{l_1,l_2,...,l_m}\right)^{\alpha}&=&\sum\limits_{l_1}p_{l_1}^{\ \alpha}\sum\limits_{l_2}p_{l_2}^{\ \alpha}...\sum\limits_{l_m}p_{l_m}^{\ \alpha} =1.
\end{eqnarray}

\section{Conclusion}

In this work, we have studied the composition of fractals. The composition is performed by alternating the iterations of the component fractals in a certain manner.
We have studied composition with several possible alternations and calculated the dimension of the composite fractals.
\begin{itemize}
\item In the case of the composition of different fractals, the composite dimension can be expressed in the form of barycenter or average of component dimensions.
\item For multifractals, the composite dimension can be expressed in a generalized Moran equation which can be solved numerically to find the composite dimension.
\item It is proved that, in general, the composite dimension is bounded by the minimum and maximum of the component dimensions, and independent of the order of different substages in a stage.
\end{itemize}
We would like to mention that the composite dimensions found in this work are true for fractals which do not have overlapping segments or where there is no elements which intersect themselves during iteration. If there is overlapping, the composite dimensions presented here constitute upper bounds for the real dimension of overlapping figures. The difference between the real dimension and the upper bound depends on the degree of intersection of the shape.

In this work we have focused on the composition by simple and regular alternation of the iterations of the component fractals. The same method can be extended to more complicated cases where the number of substages changes with the increasing iteration or where the occurrence of substages is probabilist \cite{Barnsley1}. This work is in progress. We hope that the approach of composition (or decomposition) of fractals can find application in the study of fractal geometry as well as in the modeling of different structures and processes in nature.

\end{document}